# A meshless method for the solution of incompressible flow equations

# Une méthode de résolution des problèmes d'écoulements incompressibles en absence de maillage

**Stéphane Couturier – Hamou Sadat**

*Laboratoire d'Études Thermiques, UMR CNRS 6608*
*École Supérieure d'Ingénieurs de Poitiers*
*E-mail : sadat@let.ensma.fr*

*ABSTRACT. This article reports on the efficiency of a co-located diffuse approximation method coupled with a projection algorithm for the solution of two and three-dimensional incompressible flow equations. Three typical examples show the accuracy of this meshless method.*

*RÉSUMÉ. Cet article propose de montrer l'efficacité d'une méthode de collocation basée sur l'approximation diffuse et couplée à un algorithme de projection pour la résolution des équations des écoulements incompressibles dans des domaines bi- et tri-dimensionnels. Trois exemples classiques montrent la précision de cette méthode sans maillage.*

*MOTS CLÉS : approximation diffuse, méthode sans maillage, écoulement incompressible, algorithme de projection.*

**Nomenclature**

| | |
|---|---|
| D | studied domain |
| I | functional |
| N | number of mesh points |
| P | line vector of monomials |
| $P^T$ | P-transpose |
| p | dimensionless pressure |
| Re | Reynolds number |
| t | dimensionless time |
| $(\Delta t)$ | dimensionless time step |
| u,v,w | dimensionless velocity components |
| $\vec{V}$ | dimensionless velocity vector |
| X | current point |
| x,y,z | dimensionless Cartesian coordinates |

Greek symbols

| | |
|---|---|
| $\alpha$ | vector of estimated derivatives |
| $\alpha^T$ | α-transpose |
| $\Phi$ | scalar field |
| $\nu$ | kinematic viscosity |
| $\sigma$ | practical aperture of the weighting function |
| $\omega$ | weighting function |
| $\psi$ | dimensionless stream function |

## 1. Introduction

Finite element methods [TAY 81-COM 82] and control-volume based finite element methods [PRA 85- MAS 94] are generally used in the numerical solution of fluid flow problems in regions with arbitrary shaped boundaries. For problems in which the position of large solution gradients is known a priori, such as those involving boundary layers, localized grid refinement can be used. In many situations, however, adaptive procedures for finite element meshes are necessary. This is the case in problems with moving boundaries or discontinuities for example. These procedures are time consuming and introduce numerous difficulties associated with the large number of remeshings. In recent years, in the field of computational mechanics, a new class of methods have been developed which do not require a finite element mesh. These meshless methods only require sets of discretization nodes which could be generated by several techniques such as random shooting or octree-based methods. It is worth mentioning here that CAD data structures can also be used. Nayrolles et al. [NAY 91] have presented a method called the diffuse approximation which was then used by the same authors to develop the diffuse elements method. Belytschko et al. [BEL 94] have proposed the element-free Galerkin method. An overview of these meshless methods can be found in the paper by Belytschko et al. [BEL 96]. These methods remain however relatively unpopular in the field of computational fluid dynamics. Sadat et al. [SAD 96] were the first to use the diffuse approximation to solve two-dimensional fluid flow problems. They used the vorticity-streamfunction formulation and found that the method is quite accurate and behaves well. However, streamfunction-vorticity methods are not readily extended to three dimensions. The purpose of this article is therefore to detail a primitive-variable diffuse approximation method for the calculation of incompressible fluid flow. As the enforcement of the conservation of mass is the primary challenge for CFD algorithms, we first present in the next section the projection algorithm that we used. The diffuse approximation based collocation method is then described. Three numerical examples are finally considered to show the accurcy of the method in 2D and 3D domains.

## 2. Solution procedure

Projection methods are now widely used for obtaining a steady-state solution to the momentum and continuity equations [DAT 96-ABD 87]. These procedures have proven to work well in finite element formulations. In this work, we used an equal order projection algorithm [COM 82] whose basic methodology can be summarized as follows.

The momentum and continuity nondimensional equations that govern the laminar flow of a constant propery fluid in a two- and three-dimensional domain are:

$$\frac{\partial \vec{V}}{\partial t} + \vec{V}.\nabla\vec{V} = \frac{1}{Re}\nabla^2\vec{V} - \nabla p \qquad [1]$$

$$\nabla \vec{V} = 0 \qquad [2]$$

The pressure gradient is written as sums of estimated (*) and correction (') values:

$$\nabla p = \nabla p* + \nabla p' = \nabla p_n + \nabla p' \qquad [3]$$

where the subscript (n) indicates known values from the previous step of calculation n.

These two components are associated with the two corresponding components of the velocity vector:

$$\vec{V} = \vec{V}* + \vec{V}' \qquad [4]$$

At each new time step (n+1), we write equation [1] as:

$$\frac{\vec{V}* + \vec{V}' - \vec{V}_n}{(\Delta t)} \cong \frac{1}{Re}\nabla^2\vec{V}* - \vec{V}_n.\nabla\vec{V}* - \nabla p_n - \nabla p' \qquad [5]$$

We first solve Eq. [5] for $\vec{V}*$ by writing:

$$\frac{\vec{V}* - \vec{V}_n}{(\Delta t)} = \frac{1}{Re}\nabla^2\vec{V}* - \vec{V}_n.\nabla\vec{V}* - \nabla p_n \qquad [6]$$

Since $\vec{V}*$ does not respect continuity, we can find the pressure correction p' that project $\vec{V}*$ onto the divergence free space $\vec{V}^{n+1}$ by writing:

$$\frac{\vec{V}'}{(\Delta t)} = -\nabla p' \qquad [7]$$

Taking the divergence of Eq. [7], we obtain:

$$\nabla^2 p' = -\frac{\nabla \vec{V}'}{(\Delta t)} \qquad [8]$$

where :

$$\frac{\nabla \vec{V}'}{(\Delta t)} = \frac{\nabla \vec{V}^{n+1} - \nabla \vec{V}*}{(\Delta t)} \qquad [9]$$

At this stage, continuity is enforced by writing:

$$\nabla \vec{V}^{n+1} = 0 \qquad [10]$$

We thus obtain the pressure correction equation:

$$\nabla^2 p' = \frac{\nabla \vec{V}*}{(\Delta t)} \qquad [11]$$

Eq. [11] is the Poisson equation for the pressure correction whose boundary conditions are:

$$\frac{d p'}{d n} = 0 \qquad [12]$$

at a wall boundary and

$$p' = 0 \qquad [13]$$

at a boundary where the pressure is known.

Finally, we use p' to compute the pressure field; the velocity correction obtained from Eq. [7] is then used to compute the new velocity field:

$$p^{n+1} = p* + p' \qquad [14]$$

$$\vec{V}^{n+1} = \vec{V}* + \vec{V}' \qquad [15]$$

This solution algorithm may be summarized as follows (for the calculation of the steady state solution by a false transient procedure):

1. Initial pressure and velocity fields are defined.
2. The momentum are solved to update the provisonal velocity field, Eq. [6].
3. The pressure correction equation is solved, Eq. [11].
4. The veocity corrections are calculated, Eq. [7].
5. The pressure and velocity fields are updated, Eqs. [14]-[15].
6. Increment time step and repeat steps 2-5 until convergence.

If the really transient solution is sought for, the algorithm is modified as follows :

6. Steps 2-5 are repeated until convergence
7. Increment time step and return to step 2

A convergence criterion of 0.01% or less change in all nodal values has been selected to check convergence.

The mass conservation was also monitored. At step 6, relaxations factors could be used for each variable to avoid oscillatory behaviors for high Reynolds numbers.

In a 2D domain, The flow pattern may be conveniently visualized by computing the streamlines from the equation:

$$\frac{\partial^2 \psi}{\partial x^2} + \frac{\partial^2 \psi}{\partial z^2} = \frac{\partial u}{\partial z} - \frac{\partial w}{\partial x} \qquad [16]$$

with appropriate boundary conditions.

This Poisson equation is solved a posteriori once the solution for the velocity field has been obtained.

### 3. Diffuse approximation based collocation method

Let us consider a scalar field $\Phi : R^n \to R$ whose values $\Phi_i$ are known at the points $X_i$ of a given set of N nodes in the studied domain $D \in R^n$. The diffuse approximation gives estimates of $\Phi$ and it's derivatives up to the order k at any point $X \in D$.

The key idea is to estimate the Taylor expansion of $\Phi$ at X by a weighted least squares method which uses only the values of $\Phi$ at some points $X_i$ situated in the vicinity of X.

We thus write:

$$\Phi_i^{estimated} = P^T(X_i - X)\alpha(X) \qquad [17]$$

where $P(X_i-X)$ is a vector of polynomial basis functions and $\alpha(X)$ a vector of coefficients which are determined by minimizing the quantity:

$$I(\alpha) = \sum_{i=1}^{N} \omega(X, X_i - X)\left[\Phi_i - P^T(X_i - X)\alpha(X)\right]^2 \qquad [18]$$

in which $w : R^n \to R^+$ is a rapidly decaying weight-function of compact support.

Minimization of [18] then gives:

$$A(X)\alpha(X) = B(X) \qquad [19]$$

where:

$$A(X) = \sum_{i=1}^{N} \omega(X, X_i - X)\, P(X_i - X)\, P^T(X_i - X) \qquad [20]$$

$$B(X) = \sum_{i=1}^{N} \omega(X, X_i - X)\, P(X_i - X)\Phi_i \qquad [21]$$

By inverting system [19], one obtains the components of $\alpha$ which are the derivatives of $\Phi$ at X in terms of the neighbouring nodal values $\Phi_i$.

In this work, the Taylor expansion is troncated at order k=2. We thus have:

$$\langle P(X_i - X)\rangle = \langle\, 1, (x_i - x), (y_i - y), (x_i - x)^2, \\ (x_i - x)\cdot(y_i - y), (y_i - y)^2\,\rangle \qquad [22]$$

$$\langle\alpha\rangle^T = \langle \Phi, \frac{\partial\Phi}{\partial x}, \frac{\partial\Phi}{\partial y}, \frac{1}{2!}\frac{\partial^2\Phi}{\partial x^2}, \frac{\partial\Phi^2}{\partial x \partial y}, \frac{1}{2!}\frac{\partial^2\Phi}{\partial y^2}\rangle^T \qquad [23]$$

for the two-dimensional case and:

$$\langle P(X_i - X) \rangle = \langle 1, (x_i - x), (y_i - y), (x_i - x)^2, (x_i - x) \cdot (y_i - y), \\ (y_i - y)^2, (z_i - z), (x_i - x) \cdot (z_i - z), (y_i - y) \cdot (z_i - z), (z_i - z)^2 \rangle \quad [24]$$

$$\langle \alpha \rangle^T = \langle \Phi, \frac{\partial \Phi}{\partial x}, \frac{\partial \Phi}{\partial y}, \frac{1}{2!}\frac{\partial^2 \Phi}{\partial x^2}, \frac{\partial \Phi^2}{\partial x \partial y}, \frac{1}{2!}\frac{\partial^2 \Phi}{\partial y^2}, \frac{\partial \Phi}{\partial z}, \frac{\partial \Phi^2}{\partial x \partial z}, \frac{\partial \Phi^2}{\partial y \partial z}, \frac{1}{2!}\frac{\partial^2 \Phi}{\partial z^2} \rangle^T \quad [25]$$

for the three dimensional case.

The square and symetric matrix A(X) is not singular if the number of the connected nodes at a given point is at least equal to 6 (in 2D situations) or 10 (in 3D situations). In this work, we used 9 and 27 points respectively.

Although we tried several weight-functions, the results presented in this work were obtained with the following gaussian window :

$$\omega(X, X_i - X) = \operatorname{Exp}\left(-3 \cdot \ln(10) \cdot \left(\frac{|X_i - X|}{\sigma}\right)^2\right) \quad [26]$$

$$\omega(X, X_i - X) = 0 \quad \text{if } (X_i - X)^2 > \sigma^2 \quad [27]$$

where the aperture $\sigma$ is updated at each point in order to use the same number of neighbours in the approximation.

We now turn to the solution of the partial derivatives equations of the previous section. We simply use here a collocation method. At each point of the discretization, the derivates appearing in the equation to be solved are replaced by their diffuse approximation thus leading to an algebraic system which is solved by the preconditioned biconjugate gradient method.

The Dirichlet type boundary conditions are introduced in the same way as in the finite element method. The Neumann boundary conditions on the other hand are replaced by their diffuse approximation and then introduced in the algebraic system.

**4. Numerical results**

This section is devoted to the presentation of the numerical results obtained in three test problems, namely the two- and three-dimensional lid driven cavity and the flow in a two-dimensional channel with a rectangular obstacle. The computations were performed on a personal computer (PC) with an Intel Pentium II 300-MHz processor and 64 MB of RAM.

*4.1. Two-dimensional lid-driven cavity*

The first validation case is the well-known two-dimensional lid-driven cavity flow [GHI 82]. The presented results are obtained with two values of the Reynolds number Re based on the height of the cavity and on the speed of the moving wall, namely Re=1,000 and Re=10,000

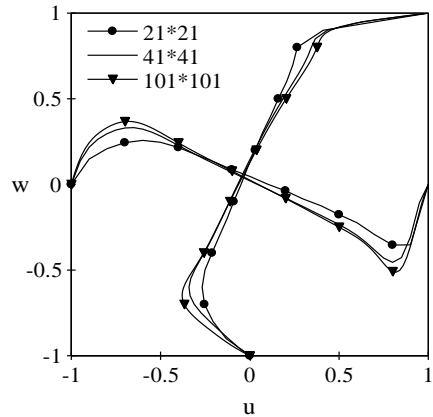

**Figure 1**. *2D lid-driven cavity : u and w velocity components along the vertical and horizontal midlines, Re=1,000.*

In order to show the spatial convergence of the method, we first show on figure 1, the velocity componenents obtained for Re=1,000 with different N×N uniform grids. The others simulations were done with a N×N non-uniform grid with a concentration of nodes near the walls. The grid is obtained by the transformation:

$$\begin{Bmatrix} x \\ y \end{Bmatrix} = \frac{1}{2} \times \left( 1 + \frac{\tanh\begin{Bmatrix} 2x_r - 1 \\ 2y_r - 1 \end{Bmatrix}}{\tanh(1)} \right) \qquad [28]$$

where $(x_r, y_r)$ correspond to the uniform grid N×N. The velocities along the midlines of the cavity are shown on figure 2. It is seen that the results are in good agreement with those of Ghia et al. [GHI 82].

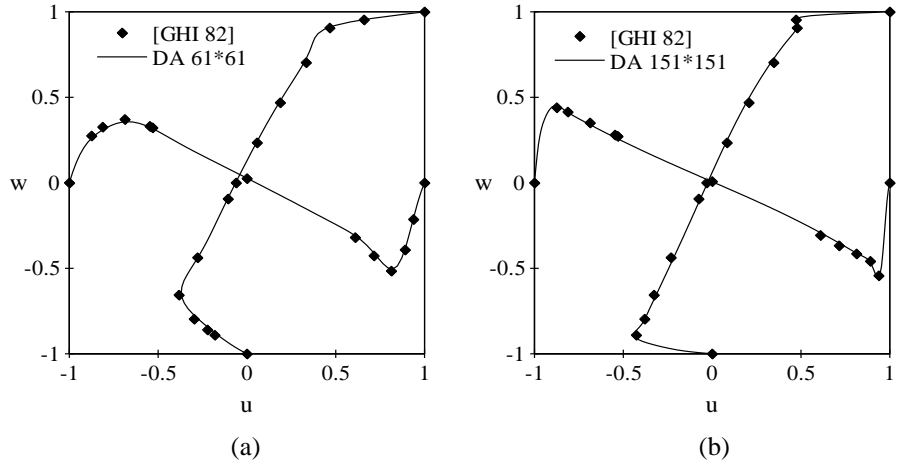

**Figure 2.** *2D lid-driven cavity. (a) Re=1,000 ; (b) Re=10,000.*

CPU time requirement and the number of iterations on various grids are reported in Table 1 for Re=1,000 (with a time step Δt=0.05 and an under-relaxation factor of 0.4 for each variable).

Figure 3 shows the streamlines and the pressure distribution. They agree with the results of [NON 97-GHI 82] very well. The values of the streamfunction at the center of the vortices obtained with a 151×151 grid are reported in Table 2. They agree very well with previously published results.

| grid size  | 31×31 | 61×61 | 101×101 | 151×151 |
|------------|-------|-------|---------|---------|
| Iterations | 1150  | 1129  | 1125    | 1126    |
| CPU (sec)  | 36    | 281   | 1571    | 6859    |

**Table 1.** *Performance comparison in terms of the number of iterations and CPU time requirement for Re=1,000.*

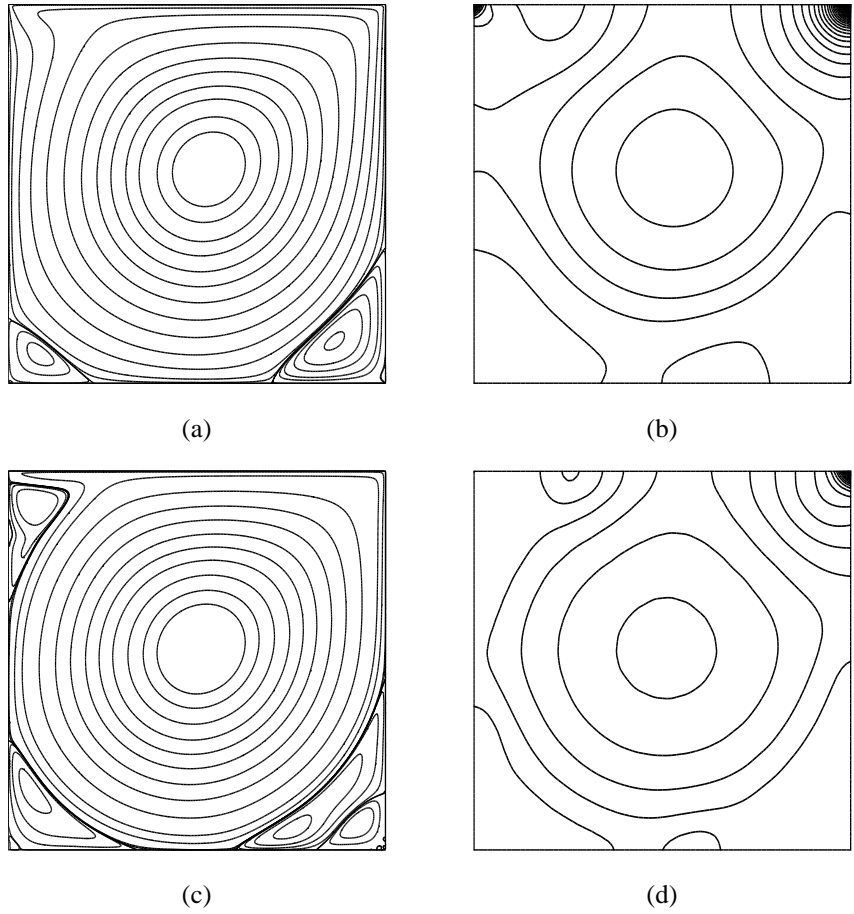

(a)    (b)

(c)    (d)

**Figure 3.** *2D lid-driven cavity. Re=1,000 : (a) streamlines, (b) pressure contour (values from -0.5 to 1.4, step 0.025). Re=10,000 : (c) streamlines, (d) pressure contour (values from -0.1 to 0.6, step 0.025).*

| Re | Primary vortex | Top left | Bottom left 1 | Bottom left 2 | Bottom right 1 | Bottom right 2 | Reference |
|---|---|---|---|---|---|---|---|
|  | -0.1184 | --- | 2.30E-4 | -3.52E-8 | 1.73E-3 | -2.50E-6 | DA |
|  | -0.1175 | --- | 2.22E-4 | -4.33E-8 | 1.72E-3 | -1.75E-7 | [NON 97] |
| 1,000 | -0.1179 | --- | 2.31E-4 | -1.14E-9 | 1.75E-3 | -9.32E-8 | [GHI 82] |
|  | -0.1151 | --- | 2.17E-4 | --- | 1.63E-3 | -1.20E-8 | [SOH 94] |
|  | -0.1190 | --- | 2.41E-4 | --- | 1.76E-3 | -6.5E-8 | [COM 94] |
|  | -0.1175 | --- | 2.32E-4 | --- | 1.77E-3 | --- | [NOB 96] |
|  | -0.1162 | 2.28E-3 | 1.47E-3 | -7.60E-7 | 3.28E-3 | -1.79E-4 | DA |
|        | -0.1199 | 2.54E-3 | 1.58E-3 | -9.55E-7  | 3.15E-3 | -1.34E-4 | [NON 97] |
| 10,000 | -0.1197 | 2.42E-3 | 1.52E-3 | -7.76E-7  | 3.42E-3 | -1.31E-4 | [GHI 82] |
|        | -0.1121 | 2.18E-3 | 1.37E-3 | -4.07E-7  | 2.80E-3 | -6.81E-5 | [SOH 94] |
|        | -0.1200 | 1.44E-3 | 1.44E-3 | -4.48E-7  | 3.34E-3 | -1.54E-4 | [COM 94] |

**Table 2**. *2D lid-driven cavity : Values of the stream function
at the center of different vortices inside the cavity.*

## *4.2. Three-dimensional lid-driven cavity*

The flow in a three-dimensional cubic lid-driven cavity is considered in this section for a Reynolds number Re=1,000.

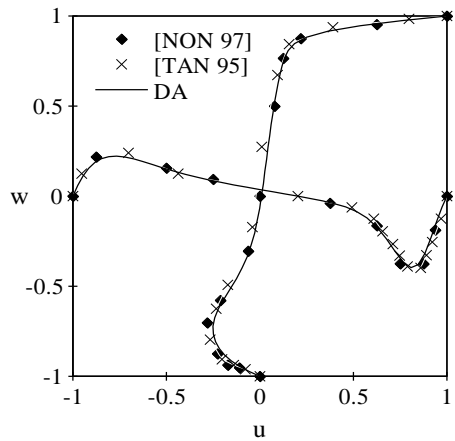

**Figure 4.** *3D lid-driven cavity at Re=1,000 : u and w velocity components
along the vertical and horizontal midlines, on the symmetry plane.*

Due to the symetry of the problem, simulations were conducted on the half of the cavity with a 41×41×21 non-uniform grid obtained by the same transformation as [28].

A meshless method for the solution of incompressible flow equations    13

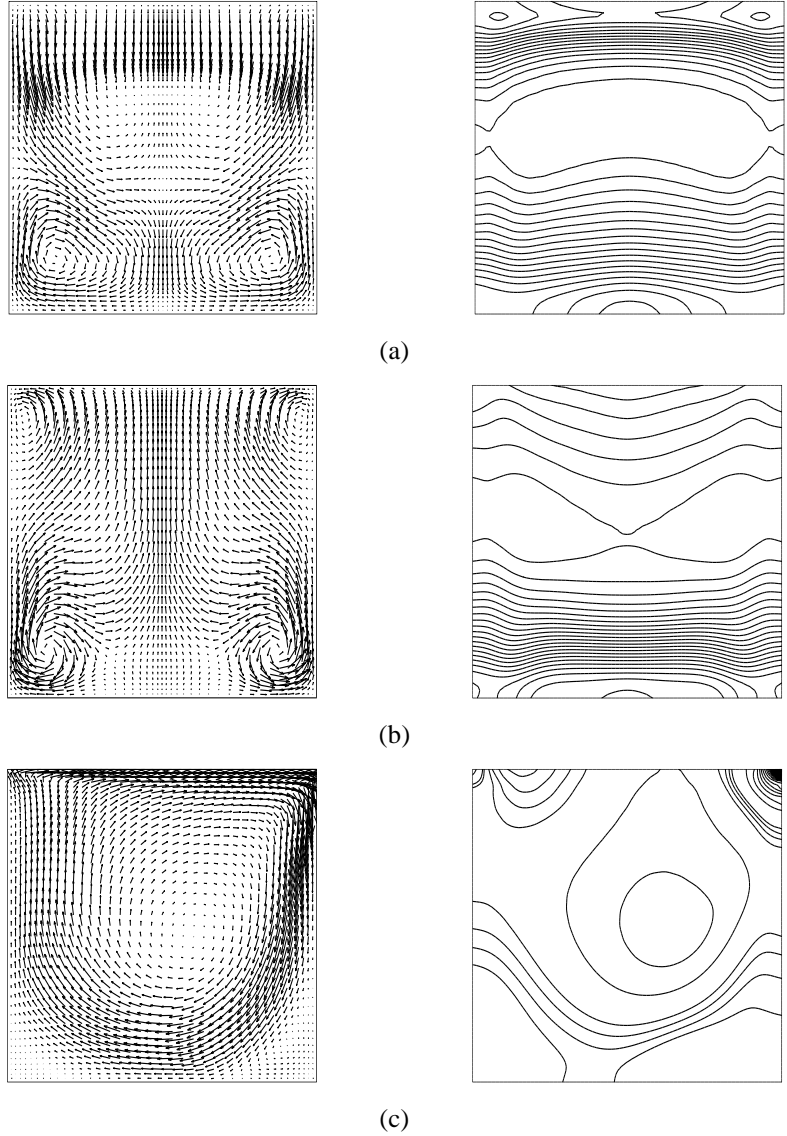

**Figure 5**. *3D lid-driven cavity at Re=1,000 : velocity vectors and pressure contours; (a) x-y midplane, (b) y-z midplane, (c) x-z midplane.*

The velocity profiles along the midlines shown on Figure 4 are in good agreement with those obtained previously [NON 97-TAN 95]. Velocity vectors and the pressure distribution in the three orthogonal midplanes are shown on Figure 5.

The x-y midplane, parallel to the moving wall is represented on Figure 5a while Figures 5b and 5c are related to the y-z transverse midplane and the x-z vertical symetry plane respectively. These results agree also satisfactorily with the literature data.

### *4.3 Flow in a two-dimensional channel with a rectangular obstacle*

This part considers the calculation of laminar flow over a rectangular block in a channel (Figure 6). The Reynolds number is defined as:

$$\mathrm{Re} = \frac{\mathrm{U}_{max}(\mathrm{H}-\mathrm{h})}{\nu} \qquad [29]$$

where (H-h) is the reference lenght and where the maximum inlet velocity $\mathrm{U}_{max}$ is the reference velocity.

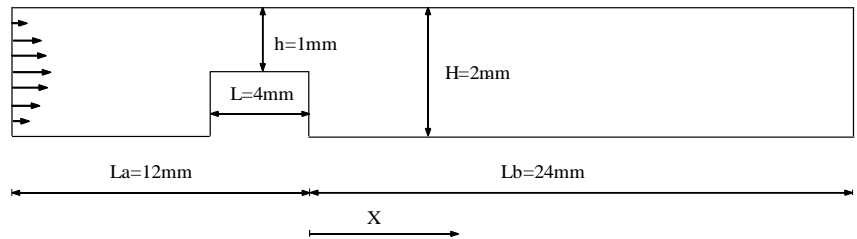

**Figure 6.** *Two-dimensional channel.*

At the inlet channel a parabolic profile for u velocity is assumed and the velocities are set to zero at the walls. At the outlet, a zero gradient is applied for the velocity while the pressure is zero.

The results obtained on a uniformly spaced 361×21 grid are compared to the numerical results of [MEL 93] and to the experimental results of [TRO 85]. Figure 7 compares predicted and measured u velocity profiles. For all positions, there is a good agreement between the results. The reattachment point is reported on Table 3. Again, a good agreement with published results is observed.

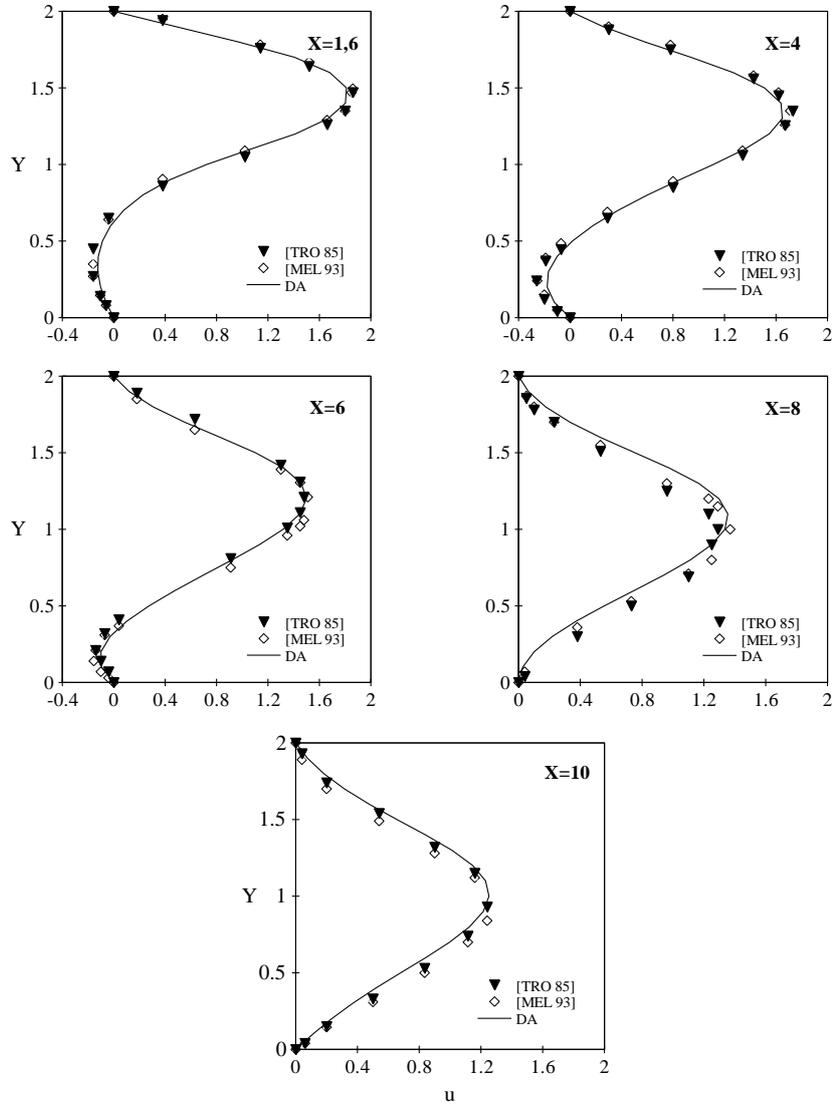

**Figure 7.** *Channel flow : u velocity profiles at different locations*.

|     | [MEL 93] | [TRO 85] | DA   |
| --- | -------- | -------- | ---- |
| $X_r$ | 7.9      | 7.1      | 7.55 |

**Table3**. *Reattachment Length*.

## 5. Conclusion

In the present work, the diffuse approximation method is presented and applied to the solution of fluid flow problems. This method provides solutions comparable in accuracy to standard numerical methods. Comparative results of test cases show good agreement and validate the applicability of the method. However, the work which has been reported is still exploratory and futher effort is needed to fully explore the limitations of the formulation.

## 6. Bibliographie